\documentclass[11pt]{article}
\usepackage{amsmath}
\usepackage{fullpage}
\usepackage{amsfonts}

\usepackage{amssymb}
\usepackage{amsthm}

\newcommand{\setlinespacing}[1]%
           {\setlength{\baselineskip}{#1 \defbaselineskip}}

\newcommand{\stan}{^{\circ}\!}

\newcommand{\dt}{\vartriangle \!\! t}

\newcommand{\nsn}{ {}^{*}\!\mathbb{N}}
\newcommand{\nsa}{ ^{*}\!A}

\theoremstyle{plain}
\newtheorem{thm}{Theorem}[section]

\newtheorem{lem}[thm]{Lemma}
\newtheorem{prop}[thm]{Proposition}
\newtheorem{defn}{Definition}[section]
\numberwithin{equation}{subsection}

\begin{document}

\title{Arithmetic progressions in Salem-type subsets of the integers}
\author{Paul Potgieter}
\date{}
\maketitle \vspace{-1.2cm}
\begin{center}
 \emph{Department of Decision Sciences, University
of South Africa\\ P.O. Box 392, Pretoria 0003, South Africa}\\
\texttt{potgip@unisa.ac.za}
\end{center}
\setcounter{section}{0}

\begin{abstract}
Given a subset of the integers of zero density, we define the weaker
notion of the fractional density of such a set. It is shown how this
notion corresponds to that of the Hausdorff dimension of a compact
subset of the reals. We then show that a version of a theorem of
{\L}aba and Pramanik on 3-term arithmetic progressions in subsets of
the unit interval also holds for subsets of the integers with
fractional density which also satisfy certain Fourier decay
conditions.


\end{abstract}

\section{Introduction}
The existence of 3-term arithmetic progressions in certain sets of
fractional Hausdorff dimension was recently established by {\L}aba
and Pramanik~\cite{LP}. They introduce Salem-type sets in $[0,1]$,
that is, sets which have a positive Hausdorff dimension and a
sufficiently rapid decay for the Fourier transform of some measure
on the set. The origins of this theorem can be traced back to Roth's
original theorem establishing 3-term arithmetic progressions in
dense subsets of the integers~\cite{Roth}. For cases where the
density of the subset is not positive, the conclusion of Roth's
theorem may still hold, providing the sets are ``random enough",
such as is the case with the primes~\cite{TaoGreen}. We will
appropriate the term ``Salem-type" to indicate a subset of the
integers which satisfy a weak density condition as well a certain
decay condition on the Fourier coefficients of its characteristic
function, as specified in Theorem 4.1.

The goal of this paper is to establish a result corresponding to
that of {\L}aba and Pramanik on the integers. The first step is to
formulate a version of Hausdorff dimension for sets which have zero
density in the conventional sense. This allows us to relax the
uniformity conditions on sets of density zero, such as discussed
in~\cite{TaoVu}.

In the second section we discuss the results that inspired this
paper. This involves a correspondence between certain subsets of
$\mathbb{N}$ and subsets of $[0,1]$. These were originally explored
by Leth~\cite{Leth}. In~\cite{Potgieter}, a nonstandard counting
formulation of Hausdorff dimension is established. Since this
formulation, when considered in the context of subsets of the
natural numbers instead of compact subsets of $\mathbb{R}$,
resembles the usual definition of density very closely, it seemed
likely that a weaker idea of density would prove useful in studying
arithmetic progressions, especially in the light of~\cite{LP}.
Indeed, when subsets of $\mathbb{N}$ are mapped to subsets of
$[0,1]$ via a mapping similar to that in~\cite{Leth}, this
``fractional density" is preserved as Hausdorff dimension.
Similarly, when a subset of $[0,1]$ is mapped into $\mathbb{N}$,
Hausdorff dimension is preserved in the guise of fractional density.

The third section discusses a uniformity condition (see for instance
\cite{TaoVu}, p161) necessary for a set of fractional density to
contain a 3-term arithmetic progression. In the fourth section, a
version of Laba and Pramanik's result is proved for subsets of
$\mathbb{N}$. The proof involves little else but repeated use of
Varnavides's theorem, as found in \cite{TaoVu}. In the final section
we construct an example of a set in the integers, analogous to that
found in Section 6 of \cite{LP}, which satisfies the conditions of
Theorem 4.1 of this paper.

Some background in nonstandard analysis is required for the second
section of this paper. A succinct but sufficient introduction to all
the necessary concepts can be found in~\cite{Potgieter}. Apart from
Definition 2.2, the rest of the paper can be read independently of
this section. However, in order to understand the motivation behind
the formulation and the direction of future investigations, it would
benefit the reader to at least give it a cursory glance.

The author would like to thank Claudius du Plooy for his generous
support during the conception of this paper, as well as for many
enlightening conversations throughout the years.

\section{Correspondence between subsets of $\mathbb{N}$ and $[0,1]$}

We use the notation of~\cite{Potgieter} throughout. Let
$A=(a_n)_{n\in \mathbb{N}}$ denote a sequence of natural numbers
(which we assume to be strictly increasing). Note that in the
paper~\cite{Leth} we are not restricted to sequences in
$\mathbb{N}$, but it will suffice for our purposes. The essential
idea behind the correspondence is to use a hyperfinite number to
divide every member of the nonstandard extension $\nsa$ of the
sequence $A$ (throughout this section we denote nonstandard
extensions of sets similarly). The standard part of a nonstandard
number or set $x$ shall be denoted by $\textrm{st}(x)$. General
results in~\cite{Leth} hold for division by any hyperfinite number
$z\in \nsn \setminus \mathbb{N}$. We shall however only consider
division of each element of $\nsa$ by the number $\langle a_n
\rangle_{\mathcal{U}}$, that is, the unique hyperfinite number
determined by the sequence $(a_n)$ under the equivalence relation of
a certain (fixed) free ultrafilter $\mathcal{U}$ (the choice of
ultrafilter is immaterial to the results). We formalise this
previous by defining:

\begin{defn} Suppose $A=(a_n)_{n\geq 1}$ is an increasing sequence
of natural numbers. Then we denote by $\textrm{st}_{z}(A)$ the set
\[ \{ \textrm{st}(a /z): a\in\, \nsa\}\] where $z=\langle a_n
\rangle_{\mathcal{U}}$.
\end{defn}

(The definition can of course be extended from $\mathbb{N}$ to
$\mathbb{Z}$.) It is clear that $st_z (A) \subseteq [0,1]$.
Furthermore, it is a closed set, as shown in Proposition 2.2
of~\cite{Leth}. Our purpose is now to show that the Hausdorff
dimension of $st_z (A)$ coincides with the ``fractional density" of
$A$. Throughout the paper we will use $[A,B]$ to denote the interval
in $\mathbb{Z}$ given by the set $\{A, A+1, \dots,B\}$. The
intervals $[A,B)$ and $(A,B]$ are defined analogously.

\begin{defn}
We say that a set $A\subseteq \mathbb{N}$ has {\emph{fractional
upper density}} $\alpha$ if
\[ \limsup_{N\to \infty} \frac{|A\cap [1,N]|}{N^{\beta}}\]
is $\infty$ for any $\beta < \alpha$ and $0$ for any $\beta >
\alpha$.
\end{defn}
(The lower fractional density can be similarly defined by replacing
$\limsup$ in the above with $\liminf$. If upper and lower fractional
densities are equal, we can just speak of the fractional density.)

We can summarise this by saying that $\bar{d}_f (A) = \alpha$. We
will sometimes need to consider the fractional density relative to a
finite but arbitrarily large number; that is, we will say that
$A\subseteq [0,N)$ has upper fractional density $\beta$ relative to
$N$ if $|A|/N^{\beta}=c>0$, for arbitrarily large $N$. In sections 3
and 4, this is mostly how the concept of upper fractional density
will be utilised. Note that the limit in the above definition is the
same as the limit of $n/a_{n}^{\beta}$, which is the form we will
use in proving Proposition 2.2.

Of course, one has to verify that such a concept yields information
that the usual definition of density does not, in the same way that
Hausdorff dimension (denoted by $\textrm{dim}_H$) yields information
that Lebesgue measure does not. Firstly, it is easily verified that
any subset of $\mathbb{N}$ of positive density has fractional
density $1$. One can also verify that there exist sets which do not
have positive density but do have positive fractional density. For
example, one can create a version of the triadic Cantor set on
$\mathbb{N}$ as follows:
\begin{itemize}
\item[1.]{Let $C_0$ be the interval $(0,3^0 ]$ (in $\mathbb{N}$). We recognise only the right-hand
endpoint of the interval, leaving $C_0 = \{1\}$.}
\item[2.]{Let $C_1$ consist of the interval $[0,3^1 ]$. Remove the middle third $(1,2]$, and keep $1$ and $3$, the right-hand
endpoints of the remaining intervals. Thus, $C_1$=\{1,3\}}
\item[3.]{Similarly with the interval $(0,3^2]$, we remove the middle third intervals $(1,2]$, $(3,4]$, $(4,5]$, $(5,6]$
and $(7,8]$. Thus, $C_2=\{1,3,7,9\}$}, and so on.
\end{itemize}

This construction can be formalised thus:
\begin{eqnarray}C_0 &=& \{1\}\\
C_{i+1} &=& C_i \cup \{ 3^{i+1}+1 - c: c\in C_i \}, \qquad i\in
\mathbb{N}\\ C&=&\bigcup_{i=0}^{\infty}C_i.
\end{eqnarray}

It is trivial to show that this set has fractional density $\log
2/\log 3$ simply by counting elements at every stage, even though it
does not have positive density.

Instead of utilising the standard definition of Hausdorff dimension,
we use the following nonstandard version~\cite{Potgieter}. Note that
for some infinitesimal $\dt = 1/N$, $N\in \nsn \setminus
\mathbb{N}$, we call the set $\{0, \dt, 2\dt, \dots, 1-\dt\}$ the
{\emph{hyperfinite time line} based on} $\dt$. The function $|\cdot
|$ denotes the transferred cardinality function.

\begin{thm}
Consider a hyperfinite time line $\mathbf{T}$ based on the
infinitesimal $N^{-1}$, for a given $N\in \nsn\setminus \mathbb{N}$.
Suppose that an internal subset $A'$ of the time line is such that
$\stan( A') = A$ and for some $\alpha >0$
\begin{eqnarray}
\stan \left( \frac{|A'|}{N^{\beta}}\right) &>& 0 \textrm{ for }
\beta
< \alpha \textrm{ and }\\
\stan \left( \frac{|A'|}{N^{\beta}}\right) &=& 0 \textrm{ for }
\beta
> \alpha.
\end{eqnarray}
Then $\alpha = \textrm{dim} A$.
\end{thm}

One might be concerned that the nonstandard formulation of Hausdorff
dimension might too closely resemble Minkowksi dimension. However,
as is shown in~\cite{Potgieter}, this formulation implies the
existence of a positive measure on a set of positive Hausdorff
dimension, which is not necessarily a property of sets of positive
Minkowski dimension.

A simple argument using the transfer principle now shows that
fractional density of the set $A$ is exactly the same as the
Hausdorff dimension of the set $\textrm{st}_{z}(A)$.

\begin{prop}
Suppose that a sequence $(a_n) = A\subseteq \mathbb{N}$ has
fractional density $\alpha$. If $z=\langle a_n
\rangle_{\mathcal{U}}$, then $\textrm{st}_z (A)$ has Hausdorff
dimension $\alpha$.
\end{prop}
{\bf Proof.} If $\beta<\alpha$, the sequence $n/a_{n}^{\beta}$ will
diverge as $n\to \infty$. Hence we can assume that for all $i$ after
a certain stage, $i/a_{i}^{\beta}>1$. If we now let $a_J$ denote the
element of the nonstandard extension of the sequence determined by
the sequence itself (modulo the free ultrafilter), the property
$J/a_{J}^{\beta}>1$ will also hold by the transfer principle.
Considering the set $\{a/a_J: a\in\, \nsa\}$, we see that it is a
subset of the hyperfinite time line based on $a_J$, since each
element of $\nsa$ is still a member of $\nsn$ (by transfer).
Furthermore, Theorem 2.1 now implies that $\textrm{dim}_H
(st_{a_J}(A))>\beta$. Similarly, for any $\beta>\alpha$, we obtain
that $\textrm{dim}_H (st_{a_J}(A))<\beta$, concluding the proof.

The converse of the previous proposition can also be easily shown by
reversing the argument, i.e. that given a subset of $[0,1]$ of
Hausdorff dimension $\alpha$, we can multiply by a hyperfinite
natural number (which is not unique) to obtain a set with fractional
density $\alpha$. A more interesting question concerns the
relationship between the Fourier-dimensional properties of compact
sets in $\mathbb{R}$ and the properties of discrete Fourier
coefficients of characteristic functions of analogous subsets of
$\mathbb{Z}$. It is this relationship we are attempting to explore
by interpreting the results in \cite{LP} in the context of the whole
numbers.

\section{Fourier conditions}

The essence of the proof of Roth's theorem, as presented in
e.g.~\cite{Lyall}, is to show that the Fourier transform of the
characteristic function of a set of positive density either
satisfies certain decay conditions, or the set has increased density
in some arithmetic progression in $\mathbb{Z}$. Iterating this
argument on the assumption that the set contains no 3-term
arithmetic progressions, a density of greater than $1$ is eventually
obtained on some arithmetic progression, a contradiction.

If the set does not have positive density, we have to impose decay
conditions on the Fourier coefficients. We first determine the
uniform rate of decay necessary to guarantee such progressions when
a set has fractional density $\alpha <1$. We will say that a subset
$A$ of a finite additive group $Z$ is $\gamma$-uniform if the
Fourier coefficients of the characteristic function satisfy
$|\widehat{\chi_A}(k)|\leq \gamma$ for all $k\in Z$, $k\neq 0$. If
this $\gamma$ is small, the set is said to be linearly uniform. In
the case of a set of positive density, it is possible to find linear
uniformity conditions which guarantee the existence of progressions.
Our version of this will be to find some $\beta$ such that if the
Fourier coefficients are all smaller than $cN^{\beta}$ for some $c$,
we will be guaranteed a 3-term arithmetic progression.

 Consider $A\subset \mathbb{Z}$ such that for some
$0< \alpha<1$, $|A\cap [0,N)| \geq \delta N^{\alpha}$ for
arbitrarily large $N$. (This implies that the upper fractional
density of $A$ is $\geq \alpha$.) We will assume, without loss, that
$|A\cap [0,N)| = \delta N^{\alpha}$ for each $N$ under
consideration. As a first approximation to the 3-term arithmetic
progressions contained $A\cap [0,N)$, we count the number of
progressions modulo $N$, i.e. the number of $x,y,z\in A$ such that
\[x+y\equiv 2z\,\textrm{mod} N.\] (In this we follow Lyall's exposition of Roth's
theorem~\cite{Lyall}, and use similar notation.) The Fourier
coefficients of a function defined on the integers modulo $N$
(denoted by $\mathbb{Z}_{N}$) are defined as usual by
\[\hat{f}(k) = \frac{1}{N}\sum_{n=0}^{N-1}f(x)e^{-\frac{2\pi ikn}{N}} \]
The number of triples satisfying the congruence, if $\chi_A$ denotes
the characteristic function of $A$, is given by
\[ \mathcal{N}_0 = N^2 \sum_{n=0}^{N-1}
\widehat{\chi_A}(n)\widehat{\chi_A}(n)\widehat{\chi_A}(-2n)\]
However, a triple satisfying the congruence does not necessarily
form a true arithmetic progression in $\mathbb{Z}$, since some of
the terms might ``wrap around" the cyclic group. If we require
instead that $x,z \in M_A = A\cap [N/3, 2N/3)$, then a
$\mathbb{Z}_N$-progression does indeed form a
$\mathbb{Z}$-progression. In this case, we estimate the true triples
$\mathcal{N}$ by writing
\[\mathcal{N}\geq N^2 \sum_{n=0}^{N-1}\widehat{\chi_{M_A}}(n)
\widehat{\chi_A}(n)\widehat{\chi_{M_A}}(-2n) = \delta N^{\alpha
-1}|M_A|^2 + N^2 \sum_{n=1}^{N-1}\widehat{\chi_{M_A}}(n)
\widehat{\chi_A}(n)\widehat{\chi_{M_A}}(-2n).\]

We require that $|M_A|\geq \frac{\delta}{4} N^{\alpha}$ and
$|\widehat{\chi_A}(k)|\leq \delta^2 N^{\beta}/32$ for $k\neq 0$.
Using the Cauchy-Schwartz inequality, this gives
\begin{eqnarray*}N^2 \left| \sum_{n=1}^{N-1}\widehat{\chi_{M_A}}(n)
\widehat{\chi_A}(n)\widehat{\chi_{M_A}}(-2n)\right| &\leq& N^2
\max_{k\neq 0}|\widehat{\chi_A}(k)| \left|
\sum_{n=1}^{N-1}\widehat{\chi_{M_A}}(n)\widehat{\chi_{M_A}}(-2n)\right|\\
&\leq& N^2 \max_{k\neq 0}|\widehat{\chi_A}(k)|
\left(\sum_{n=1}^{N-1}|\widehat{\chi_{M_A}}(n)|^2
\right)^{\frac{1}{2}} \left(
\sum_{n} |\widehat{\chi_{M_A}}(-2n)|^2\right)^{\frac{1}{2}}\\
&\leq& N^2\max_{k\neq 0}|\widehat{\chi_A}(k)| \sum_{n=0}^{N-1}
|\widehat{\chi_{M_A}}(n)|^2\\
&=&N^2\max_{k\neq 0}|\widehat{\chi_A}(k)|\cdot \frac{1}{N}\sum_{x=0}^{N-1} \chi_{M_A} (x)\\
&\leq & \frac{\delta^2}{32}N^{\beta} \cdot N \cdot |M_A|\\
&\leq&\frac{\delta^3}{32}N^{\beta} N^{\alpha +1}
\end{eqnarray*}

If we now require that $\beta < 2\alpha -2$, say $\beta = 2\alpha -2
-\varepsilon$, we find that

\begin{eqnarray*}\mathcal{N}&\geq& \delta N^{\alpha -1}|M_A|^2 -
\frac{\delta^3}{32}N^{3\alpha -1 -\varepsilon}\\
&\geq& \frac{\delta^3 N^{3\alpha -1}}{32}(2-N^{-\varepsilon})\\
&>&\frac{\delta^3 N^{3\alpha -1}}{32}
\end{eqnarray*}
This will be large for $\alpha > 1/3$.

We still have not taken into account the number of trivial
progressions $x=y=z$, of which there are $|A|=\delta N^{\alpha}$. If
we subtract this from the estimate obtained above and require that
$\alpha > 1/2$ and (for instance) $N> 32/\delta^2$, we are certain
to have a non-trivial 3-progression.

Of course, we might not always be as fortunate as to have such small
Fourier coefficients. In the next section we show that weaker
non-uniform conditions would still suffice, provided that the decay
is sufficiently structured.

\section{Salem-type sets in the integers}

In this section we prove the following:

\begin{thm}
Let $A\subseteq \mathbb{Z}$. Suppose $A$ satisfies the following
conditions:
\begin{itemize}
\item[(i)]{$A$ has upper fractional density $\alpha$, where $\alpha >1/2$.}
\item[(ii)]{The Fourier coefficients of the characteristic functions $\chi_{A_N}$ of $A_N = A\cap[0,N-1]$ satisfy
\[|\widehat{\chi_{A_N}} (k)| \leq C(|k|N)^{-\beta/2}\] for large $N$, for some $2/3< \beta \leq
1$ satisfying $\beta>2-2\alpha$.}
\end{itemize}
Then $A$ contains an arithmetic progression of length 3.
\end{thm}

As long as the interval $[0,N-1]$ is fixed, as it is throughout most
of the proof, we will use simply $A$ instead of $A_N$.

To prove Proposition 4.1, we use a modified version of the density
arguments using Varnavides's theorem, to be found in
e.g.~\cite{TaoVu}. Throughout, we use $Z$ to denote a finite
additive group of odd order $N$. The expectation of a function on
$Z$ is defined as
\[ \mathbf{E}_Z (f) = \mathbf{E}_{x\in Z}(f) = \frac{1}{|Z|}\sum_{x\in
Z}f(x)\] The $L^p (Z)$-norm of a function $f:Z \to \mathbb{C}$ is
given by
\[\|f\|_{L^p (Z)} = \left( \frac{1}{N}\sum_{n=0}^{N-1}|f(n)|^p
\right)^{\frac{1}{p}}.\] We also define the {\it linear bias} of a
function $f:Z\to \mathbb{C}$ by
\[ \|f\|_{u^2 (Z)}= \sup_{\xi \in Z} |\hat{f} (\xi
)|.\]

In the proof we will repeatedly use the following definition:
\begin{defn}
\[ \Lambda_3 (f,g,h) = \mathbf{E}_{x,r\in Z} f(x)g(x+r)h(x+2r)        \]
\end{defn}
Note that $|Z|^2\Lambda_3 (\chi_A ,\chi_A,\chi_A)$ is an indication
of the number of 3-term arithmetic progressions to be found in a set
$A\subseteq Z$, although some might be counted more than once. To
remove trivial progressions, one has to subtract $|A|$. It follows
that if $|Z|^2\Lambda_3 (\chi_A ,\chi_A,\chi_A)-|A|$ is suitably
large, $A$ will contain at least one 3-progression as a subset of
the group $Z$.

The following can be found in~\cite{TaoVu}, p.374.

\begin{prop}
For functions $f$, $g$ and $h$ from $Z$ to $\mathbb{C}$,
\begin{equation} \Lambda_3 (f,g,h) =  \sum_{n=0}^{N-1}
\widehat{f}(n)\widehat{g}(-2n)\widehat{h}(n).\end{equation}
\end{prop}

We also have the following property of $\Lambda_3$~\cite{TaoVu}:
\begin{equation} \Lambda_3 (f,g,h) \leq \|f\|_{u^2 (Z)}\|g\|_{L^2
(Z)}\|h\|_{L^2 (Z)}
\end{equation}

{\bf Proof of Theorem 4.1.} We can assume that not all of the
Fourier coefficients are smaller than or equal to $\delta^2
/8N^{\beta}$, since that would immediately imply a 3-term arithmetic
progression, by the result in Section 3.

From now on, we denote $\chi_A$ by $\mu$, for brevity and also to
consolidate the analogy with~\cite{LP}. Where they consider a
compact subset of $[0,1]$ of certain Hausdorff dimension together
with a sufficient decay of the measure guaranteed to exist on the
set, we consider a set of fractional density with sufficient decay
of the discrete Fourier transform of the characteristic function.

We decompose $\mu$ into a sum $\mu_1 +\mu_2$. Using this, we
estimate the expression $\Lambda (\mu ,\mu ,\mu)$. If this is large
enough, it will guarantee the existence of a 3-term arithmetic
progression.

We let $F_K$ denote a version of the Fej\'{e}r kernel on $[0,N-1]$:
\[ F_K (x) = \sum_{k=0}^{K}\left( 1-\frac{k}{K+1} \right)e^{\frac{2\pi
ikx}{N}}.\]

Define $\mu_1$ as the convolution of $\mu$ and $F_K$:
\[ \mu_1 (x)=(F_K \ast \mu )(x) = \sum_{y=0}^{N-1} \sum_{n=0}^{K}
\left( 1-\frac{n}{K+1}\right) e^{\frac{2\pi in(x-y)}{N}} \mu(y).\]

By rewriting the convolution product, we can find the Fourier series
of $\mu_1$:

\begin{eqnarray*}
\mu_1 (x)&=& \sum_{k=0}^{K} \sum_{n=0}^{N-1} \left(
1-\frac{k}{K+1}\right) e^{\frac{2\pi ikx}{N}}e^{-\frac{2\pi ikn}{N}}\chi_A (n)\\
&=&  \sum_{k=0}^{K} \left( 1-\frac{k}{K+1} \right) e^{2\pi i kx}
\widehat{\chi_A}(k).
\end{eqnarray*}

Thus, if $n<K+1$,
\[ \widehat{\mu_1}(n)= \left(1-\frac{n}{K+1} \right)
\widehat{\chi_A}(n).\] Otherwise, $\widehat{\mu_1}(n)=0$. Also,
since $\widehat{\mu_2} (n) = \widehat{\chi_A}(n) -
\widehat{\mu_1}(n)$,
\[\widehat{\mu_2}(n) = \textrm{min}\left(1, \frac{n}{K+1}\right)
\widehat{\chi_A}(n).\]

To calculate $\Lambda_3 (\mu ,\mu ,\mu)$, we split the expression
$\Lambda_3 (\mu_1 +\mu_2 ,\mu_1 +\mu_2 ,\mu_1 +\mu_2)$ into eight
terms of the form $\Lambda_3 (\mu_i ,\mu_j ,\mu_k)$, $i,j,k \in
\{1,2\}$. The idea is then to show that the term $\Lambda_3 (\mu_1
,\mu_1,\mu_1)$ dominates the others, and will be large enough to
guarantee an arithmetic progression.

We can now use the following inequality, which follows from (4.2):
\begin{equation}|\Lambda_3 (f ,g,h)|\leq \sum_{0\leq n< N}
|\hat{f}(n)||\hat{g}(-2n)||\hat{h}(n)|.\label{eq:ineq1}\end{equation}

We only evaluate two of the terms which contain at most two
instances of $\mu_1$. The others can be evaluated according to the
exact same principles.

Throughout the calculations, we assume that $K<N/2$, so that, for
$1\leq n \leq K$, $|-2n|^{-\frac{\beta}{2}} \leq
(2n)^{-\frac{\beta}{2}}$. Furthermore, this implies that $\min\{1,
1-(N-2n)/(K+1)\}=1$ for $1\leq n\leq K$. We will later see that the
lower bound we place on $K$ does not violate these conditions. This
assumption allows us to replace $|-2k|$ by $2k$ in the sequel.

First considering the term $\Lambda_3 (\mu_1 ,\mu_2 ,\mu_1)$, we
know from inequality~\ref{eq:ineq1}, the fact that $\widehat{\mu_1}
(n) =0$ for $n\geq K+1$ and $\widehat{\mu_2}(0)=0$ that
\begin{eqnarray*}|\Lambda_3 (\mu_1 ,\mu_2 ,\mu_1)|&\leq& \sum_{0<
n\leq N}
|\widehat{\mu_1}(n)|^2|\widehat{\mu_2}(-2n)|\\
&=& O\left( N^{-\frac{3\beta}{2}}\sum_{0< n\leq
K}\left(1-\frac{n}{K+1}\right)
n^{-\frac{3\beta}{2}} \right)\\
&=&O\left( N^{-\frac{3\beta}{2}}\sum_{0< n\leq K}
n^{-\frac{3\beta}{2}} \right)\\
&=& O\left( N^{-\frac{3\beta}{2}} \right),
\end{eqnarray*}
since the sum $\sum_{n=1}^{\infty} n^{-\frac{3\beta}{2}}$ is
convergent.

Next, we turn to the expression $\Lambda_3 (\mu_1 , \mu_2 ,\mu_2)$.
Using the same properties of the Fourier coefficients, we find once
again that
\begin{eqnarray*}
|\Lambda_3 (\mu_1 , \mu_2 ,\mu_2)|&\leq&\sum_{0< n\leq N}
|\widehat{\mu_1}(n)||\widehat{\mu_2}(-2n)||\widehat{\mu_2}(n)|\\
&=&O\left( N^{-\frac{3\beta}{2}}\sum_{0< n\leq
K}\left(1-\frac{n}{K+1}\right) n^{-\frac{3\beta}{2}} \right).
\end{eqnarray*}
The same bound clearly applies as for the previous expression.
Because the Fourier coefficients of $\mu_1$ are $0$ for $n\geq K+1$,
any term involving $\mu_1$ can be approximated this way. If the term
does not involve $\mu_1$, we have no such cut-off, yet even without
such we can still easily obtain an upper bound of
$O(N^{-\frac{3\beta}{2}})$ on $|\Lambda_3 (\mu_2,\mu_2,\mu_2)|$.

Hence, all terms in the expansion of $\Lambda_3
(\mu_1+\mu_2,\mu_1+\mu_2,\mu_1+\mu_2)$ that involve $\mu_2$ become
at most $O(N^{-\frac{3\beta}{2}})$. The next step is to show that
the term $\Lambda_3 (\mu_1,\mu_1,\mu_1)$ is large compared to these.

To do so, we once again decompose the relevant function into two
parts. Set
\begin{eqnarray*} \mu_3 &=& \mu_1 -\mathbb{E}(\mu_1)\\ \textrm{and }
\mu_4 &=&\mathbb{E}(\mu_1).
\end{eqnarray*}
We approximate the expression
\[\Lambda_3 (\mu_1,\mu_1,\mu_1) = \Lambda_3 (\mu_3+\mu_4,
\mu_3+\mu_4, \mu_3+\mu_4)\] by showing that one term is large
compared to the seven others.

It is clear that $\Lambda_3 (\mu_4,\mu_4,\mu_4) = \delta^3
N^{3\alpha -3}$. Furthermore, $\widehat{\mu_3}(0) =0$ and
$\widehat{\mu_3}(k)=\widehat{\mu_1}(k)$ for $k>0$. As in the
previous part of the proof, we now use inequality 4.8 to approximate
the lesser terms. Firstly,
\begin{eqnarray*}
|\Lambda_3 (\mu_3,\mu_4,\mu_3)|&\leq& \|\mu_3 \|_{u^{2}(Z)}\|\mu_4
\|_{L^{2}(Z)}\|\mu_3\|_{L^2(Z)}\\
&=& \delta N^{\alpha -1}\left[ \max_n
\left\{\left|\left(1-\frac{n}{K+1}\right)\widehat{\chi_A}(n)\right|
\right\}\right] \|\mu_3\|_{L^2 (Z)}\\
&=& O(\delta N^{\alpha -\frac{\beta}{2} -1}\|\mu_3\|_{L^2 (Z)}).
\end{eqnarray*}

Assuming that $K=O(N^{\frac{1}{3}})$, we can use Parseval and  an
integral to approximate the $L^2(Z)$-norm of $\mu_3$:
\begin{eqnarray*} \|\mu_3 \|^{2}_{L^2 (Z)}&\leq& \sum_{n=0}^{N-1}
|\widehat{\mu_3}(n)|^2\\
&=& O\left( \sum_{n=1}^{K}\left( 1-\frac{n}{K+1}\right)
k^{-\beta}N^{-\beta}\right)\\
&=& N^{-\beta} O\left[ \int_{1}^{K} \left(
1-\frac{x-1}{K+1}\right)^2
(x-1)^{-\beta}dx + \left(1-\frac{1}{K+1}\right)\right]\\
&=&O(N^{-\frac{4\beta}{3}+\frac{1}{3}}).
\end{eqnarray*}

Therefore, $\|\mu_3\|_{L^2 (Z)} =
O(N^{-\frac{2\beta}{3}+\frac{1}{6}})$. It follows that the term
$|\Lambda_3 (\mu_3,\mu_4,\mu_3)|$ is
$O(N^{-\frac{19}{6}+\frac{10\alpha}{3}})$ (remembering that
$\beta>2-2\alpha$, and the same clearly holds for $|\Lambda_3
(\mu_3,\mu_3,\mu_4)|$. Similar calculations show that similar upper
bounds hold for every term involving $\mu_3$. Since $\alpha >1/2$,
this is small compared to $N^{3\alpha -3}$.

All of the approximations now imply that
\[\Lambda_3 (\mu,\mu,\mu) = \Omega(N^{3\alpha - 3}).\]

The number of arithmetic progressions in $Z$ is counted by the
expression
\[ N^2 \Lambda_3 (\chi_A, \chi_A,\chi_A) - |A|\]
where the second term is employed to ensure we disregard
progressions with difference $0$. It is important to observe here
that the progressions counted is the number of proper progressions
(i.e. with non-zero difference) in the cyclic group $Z$, which may
not be equivalent to the number of progressions in the interval
$[0,N-1]\subset \mathbb{Z}$ (which will be referred to as
\emph{genuine} progressions). The question is now how to eliminate
the progressions which ``wrap around" the cyclic group $Z$. In
Roth-type theorems, this is often done through density-increment
arguments, for instance in chapter 10 of \cite{TaoVu}. In our case,
we instead consider the set $A$ as a subset of the interval
$[0,3N)$, which we can again consider as a cyclic group, which we
will call $Z'$. (This is an embedding of $A_N$ into $[0,3N)$, not a
restriction of the original set to a larger interval.) Any proper
progression in $A$, seen as a subset of $Z'$, would now have to be a
genuine progression, since there are no elements of $A$ in the
interval $[N,3N)$. Assuming that there are no progressions except
trivial ones, this means that the total number is simply the
cardinality of $A$.

We still denote the characteristic function of $A$ as a subset of
$Z$ by $\chi_A$, whereas the characteristic function of $A$ as a
subset of $Z'$ is denoted by $\chi_{A'}$. The effect on the Fourier
coefficients of $\chi_{A}$ is to ``smear" them in such a way that
their contribution to the sum-of-squares in the Parseval inequality
is taken up by several Fourier coefficients of $\chi_{A'}$. By
simply using the definition of the Fourier coefficients, it is
easily shown that
\[|\widehat{\chi_{A'}}(3k)|^2+|\widehat{\chi_{A'}}(3k-1)|^2+|\widehat{\chi_{A'}}(3k-2)|^2
\leq \frac{1}{3}|\widehat{\chi_{A}}(k)|^2.\] This now has the
implication that $A'$ satisfies condition (ii) of Theorem 4.1, with
some slightly modified constants. It is also obvious that $A'$ has
the same fractional density $\alpha$ as $A$. Thus, the proof implies
that the number of three-term arithmetic progressions in $A'$ is
greater than
\[cN^{3\alpha-1}-|A'|\] for some constant $c$. Since all
progressions counted by this expression are genuine, we have
established the existence of the required progressions in $A$.

\section{Example of a Salem-type set}

In this section we present a version in the whole numbers of the
Salem-type set constructed in~\cite{LP}.

Consider the set $\{0,1,2,\dots ,N^j-1\}$ for $N$ and $j$ large, and
some $t$, $1\leq t \leq N$. Our aim is to construct a set which has
fractional density $\alpha=\log t/\log N$ (relative to the finite
set $N^{j}$) and for which the Fourier coefficients of the
characteristic function satisfy condition (ii) of Proposition 4.1,
with $\beta > 2-2\alpha$. At each of the $j$ stages of the
construction, we randomly pick a number of points from the total in
a ratio $t/N$, in such a way that the Fourier coefficients of
successive sets satisfy certain inequalities.

Let $A_0= \{0,1,\dots, N^j-1\}$. Divide $A_0$ into $N$ equal
intervals (in the whole numbers, as usual) of length $N^{j-1}$. Let
the left-hand endpoints of these intervals be denoted by \[B^{*}_{0}
= \{0, N^{j-1},2N^{j-1}  \dots,(N-1)N^{j-1}\}.\] From this set we
choose $t$ elements with equal probability $1/t$, and call this
$B_{0}$. We form $A_1$ from this by setting
\[ A_1 = \bigcup_{b\in B_{0}}\{b, b+1, \dots, b+N^{j-1}-1\}.\]
We now divide each interval of $A_1$ into $N$ equal pieces of length
$N^{j-2}$ and form the set \[B^{*}_{1} = \bigcup_{b\in B_{0}} \{b,
b+N^{j-2}, \dots, b+(N-1)N^{j-2}\}\] from the endpoints of the
intervals newly divided. For each of the $t$ components in the union
constituting  $B^{*}_{1}$, we now have $N$ elements, and from each
choose $t$ uniformly and call the resulting (random) set $B_1$.  The
choice of $t$ elements associated to an element $b$ of $B^{*}_{1}$
we call $B_{x(b)}$, whilst the portion of $B^{*}_{1}$ of length
$N^{j-2}$ starting at $b$ is denoted by $B^{*}_{1,b}$. Iterating
this construction, we obtain from a set $A_m$ consisting of $t^m$
intervals of length $N^{j-m}$, a subdivision characterised by
$B^{*}_{m+1}$ and a choice of $t^{m+1}$ subintervals characterised
by $B_{m+1}$, which we then use to obtain $A_{m+1}$.

Some quick calculation will show that this set has fractional
density $\log t/\log N$ relative to each interval $[0,N^j)$. In
order to determine the rate of decay of the discrete Fourier
transform, we borrow the technique utilised in \cite{LP}, pp.
20--26, adapted to the whole numbers. Fundamental to the calculation
is a version of Bernstein's inequality by Ben Green~\cite{Green}.

\begin{lem}
Let $X_1, \dots ,X_n$ be independent random variables with
$|X_j|\leq 1$, $\mathbb{E}X_i =0$ and $\mathbb{E}|X_j|^2
=\sigma_{j}^{2}$. Let $\sum \sigma_{j}^{2} \leq \sigma^2$, and
assume that $\sigma^2 \geq 6n\lambda$. Then
\[ \mathbb{P}\left( \left| \sum_{1}^{n} X_j \right| \geq
n\lambda\right) \leq 4e^{-n^2\lambda^2 /8\sigma^2}.\]
\end{lem}

Given a set $B\subset [0,1]$, we write
\[ S_B (k) = \sum_{b\in B} e^{-2\pi ikb}.\]
If we are instead considering a set $B\subset \mathbb{Z}$ with
$B\subset [0,N^j)$, we abuse the notation by also using $S_B (k)$ to
denote the sum \[\sum_{b\in B} e^{-\frac{2\pi ikb}{N^j}}.\] In this
way, we can either regard $S_B$ as an exponential sum, or as the
Fourier transform of the characteristic function multiplied by a
factor $N^j$.

The previous lemma can be used to prove the following, which is a
restatement of Lemma 6.2 in \cite{LP}:

\begin{lem}
Let $B^* = \{0,\frac{1}{MN},\frac{2}{MN},\dots,\frac{N-1}{MN}\}$ and
let $1\leq t \leq N$. Let
\[ \eta^2 t=32 \log{8N^2M}\] Then there exists a set $B(x)\subset B^*$
with $|B|=t$ such that
\[ \left| \frac{S_{B(x)}(k)}{t}-\frac{S_{B^*}(k)}{N}\right| \leq
\eta \,\textrm{ for all } k \in [0,MN),\, x\in [0, N-1],\] where
\[ B(x) = \left\{ \frac{(x+y)mod N}{MN}: y\in B\right\}.\]
\end{lem}

In the proof of this from Lemma 5.1, it is shown that the condition
is satisfied with probability greater than half, indicating that at
least half of all possible choices of $B(x)$ will have the property.

One more tool will be necessary before we start the proof -- an
approximation of the Fourier coefficients by an integral.
Specifically, by considering the integral of a smooth function
$f:\mathbb{R}\to \mathbb{C}$ from $a$ to $b$ as being approximated
by a left Riemann sum with step-size $\Delta = (b-a)/M$, we get

\[\left| \int_{a}^{b} f(x)dx - \Delta \sum_{n=0}^{M-1}f(a+n\Delta)\right| \leq \frac{c(b-a)^3}{M^2}
\sup_{x\in [a,b]}|f''(x)| ,\]

where the constant $c$ is independent of $M$, $a$ and $b$.

We can now use a proof similar to that in \cite{LP}, with some
adjustment for the error term.

Define
\[ \psi_m (k)= \frac{N^{m}}{t^m}\widehat{\chi_{A_m}}(k) =
\frac{N^{m}}{t^m}\left( \frac{1}{N^j} \sum_{a\in A_m}e^{-\frac{2\pi
ika}{N^j}}\right).\] Although $\psi_m$ is not quite the same as the
Fourier transform, it will yield enough information to determine an
upper bound.

Let $B_m$ be in relation to $A_m$ as in the construction above. Then
\begin{equation} \psi_m (k) = \frac{N^{m}}{t^m} \sum_{b\in B_m}
\frac{1}{N^j}\left( e^{-\frac{2\pi ikb}{N^j}}+e^{-\frac{2\pi
ik(b+1)}{N^j}}+\cdots +e^{-\frac{2\pi
ik(b+N^{m-j}-1)}{N^j}}\right)\end{equation}

Note that if the left-hand endpoint of a subinterval of length
$N^{j-m-1}$ is determined, the whole interval is determined. If we
consider a choice of $t$ numbers from a collection of $N$ numbers to
determine the start of the interval, the exact same choice can be
considered to be applied $N^{j-m-1}$ times, from a sample space
consisting of translates of the $N$ starting points of the
intervals. In the Fourier transform of the characteristic function
of the interval, these terms then contribute the same as the
starting point, except for a phase shift for each element. If we now
wish to compute the difference $|\psi_{m+1}-\psi_m|$, the above
expression for $\psi_m$ shows that we can consider the difference

\begin{equation} \left| \frac{N^{m+1}}{t^m} \sum_{b\in
B_m}\left|\frac{S_{B^{*}_{m,b}}(k)}{N}
-\frac{S_{B_{x(b)}}(k)}{t}\right| \left(\frac{1}{N^j}
\sum_{n=0}^{N^{j-m-1}-1} e^{-\frac{2\pi ikn}{N^j}}\right)\right|
\end{equation}

In the above, we stay close to the notation of \cite{LP} in denoting
the exponential sum over the set $B^{*}_{m,b} = \{b, b+N^{j-m-1},
b+2N^{j-m-1},\dots, b+(N-1)N^{j-m-1}\}$ by $S_{B^{*}_{m,b}}$ and the
sum over the corresponding $t$-choice by $S_{B_{x(b)}}$. We now
approximate the final sum by an integral:
\[\frac{1}{N^j}\sum_{n=0}^{N^{j-m-1}-1}e^{-\frac{2\pi
ikn}{N^j}}=\int_{0}^{N^{-(m+1)}}e^{-2\pi ikx}dx + O\left(\frac{k^2
N^{-3(m+1)}}{N^{2j}}\right),\] where the error term is that of a
Riemann sum-approximation of the integral using a step-size
$N^{-j}$.

The error term can easily be shown to be less than the integral in
absolute value, especially keeping in mind that we can choose $N$
arbitrarily large. Hence we dispose of it in the absolute value,
keeping in mind that it might necessitate the use of a constant
$c<2$, which is not dependent on $m$. Computing the integral, we
find
\begin{equation} |\psi_{m+1}(k)-\psi_m (k)|\leq c\frac{(1-e^{-2\pi ik/N^{m+1}})}{t^m(2\pi
ik/N^{m+1})}\sum_{b\in B_m}\left|\frac{S_{B^{*}_{m,b}(k)}}{N}
-\frac{S_{B_{x(b)}}(k)}{t}\right|.\end{equation}

It is now obvious that the above equation is very nearly of the same
form as (52) in Lemma 6.4 of \cite{LP}. We can therefore apply the
result of the lemma to obtain
\begin{equation}|\psi_{m+1}(k)-\psi_m (k)|\leq 32 \min\left(1,
\frac{N^{m+1}}{|k|}\right)t^{-\frac{m+1}{2}}\log(8N^{m+1}).
\end{equation}

We now show that the condition 4.1 (ii) is satisfied for any $\beta
>\alpha$ such that $\beta >2-2\alpha$. Since $\psi_0 (k)=0$ for all $k\in
\{0,1,\dots, N^{j}-1\}$, we can find an upper bound on $\psi_j (k)$
by bounding the sum of all such differences. By noting that $t =
N^{\alpha}$, we can write the summand as follows (ignoring the
constant factor, which has no bearing from here on):

\begin{eqnarray}
\min\left(1, \frac{N^{m}}{k}\right)t^{-\frac{m}{2}}\log(8N^{m})&=&
\min\left(1,
\frac{N^{m}}{k}\right) N^{-\frac{\alpha m}{2}}(\log 8+m\log N )\\
&=&\min\left(1, \frac{N^{m}}{k}\right)N^{-\frac{\beta m}{2}}
N^{-\frac{(\alpha -\beta) m}{2}}(\log 8+m\log N )
\end{eqnarray}

Using the fact that $N^{-(\alpha -\beta)m/2}j\log N\leq 2(\alpha
-\beta)^{-1}$ (which can be established using elementary
calculus~\cite{LP}), the sum is bounded by

\begin{equation}
\sum_{m=1}^{j}\min \left(1,\frac{N^m}{k}\right)N^{-\frac{\beta
m}{2}}\left(N^{-\frac{(\alpha -\beta) m}{2}}\log 8+2(\alpha
-\beta)^{-1}\right) \leq \sum_{m=1}^{j}\min
\left(1,\frac{N^m}{k}\right)N^{-\frac{\beta m}{2}}\left(\log
8+2(\alpha -\beta)^{-1}\right)
\end{equation}

We consider two different regions: one where $1\leq m
\leq \log k/\log N$ and one where $m> \log k/\log N$. In the first
case,

\begin{equation}
S_1 = k^{-1}t^j(\log 8 +2(\alpha -\beta)^{-1})\sum_{1\leq
m\leq\frac{\log k}{\log N}} N^{m(1-\frac{\beta}{2})}
\end{equation}
The sum on the right is easily bounded, thus
\begin{equation}
S_1 \leq 2k^{-1}(\log 8 +2(\alpha
-\beta)^{-1})k^{1-\frac{\beta}{2}}\leq C_1 k^{-\beta/2}
\end{equation}
for some $C_1$ independent of $N$, $j$.

Approximating the second part of the sum is similar, and we obtain
\begin{equation}
S_2 =(\log 8 +2(\alpha -\beta)^{-1})\sum_{\frac{\log k}{\log
N}<m\leq j} N^{-\beta m/2}\leq C_2 k^{-\beta /2}.
\end{equation}
Using the bounds for $S_1$ and $S_2$, we get
\[|\psi_j (k)|\leq C|k|^{-\beta/2}.\]

We can obtain $\widehat{\chi_{A_j}}(k)$ by multiplication of $\psi_j
(k)$ by a factor $t^j /N^j$. Because of the construction,
\[t^j/N^j=N^{(\alpha-1)j}<N^{-\frac{\beta j}{2}},\] since we chose
$\beta>2-2\alpha$. This yields the desired bound on the Fourier
coefficients.

By this example and the result in the previous section, there seems
to be a clear correspondence between perfect subsets of $[0,1]$ and
sets in $\mathbb{Z}$, which preserves Hausdorff- and
Fourier-dimensional properties. An examination of the precision of
the correspondence will appear in the sequel to this paper.

%%%%%%%%%%%%%%%%%%%%%%%%%%%%%%%%%%%%%%%%%%%%%%%%%%%%%%%%%%%%%%%%%%%%%%%%%%%%%%

\end{document}